\documentclass[11pt]{article}
\usepackage{lmodern}
\usepackage[margin=1in]{geometry}
\usepackage{amsmath,amssymb,amsthm}
\usepackage{seqsplit}
\usepackage{enumitem}
\usepackage{microtype}
\usepackage{times}
\usepackage{xcolor}
\usepackage[colorlinks=true,linkcolor=blue!50!black,citecolor=blue!50!black,urlcolor=blue!50!black]{hyperref}

\newtheorem{theorem}{Theorem}[section]
\newtheorem{lemma}[theorem]{Lemma}
\newtheorem{proposition}[theorem]{Proposition}

\newcommand{\cC}{\mathcal C}
\newcommand{\cZ}{\mathcal Z}
\newcommand{\cO}{\mathcal O}
\newcommand{\one}{\mathbf 1}

\newcommand{\Hom}{\operatorname{Hom}}
\newcommand{\Tr}{\operatorname{Tr}}
\newcommand{\Gal}{\operatorname{Gal}}

\newcommand{\ord}{\operatorname{ord}}
\newcommand{\QQ}{\mathbb Q}
\newcommand{\CC}{\mathbb C}
\newcommand{\ZZ}{\mathbb Z}
\newcommand{\NN}{\mathbb Z_{\geq 0}}

\begin{document}
\title{The multiplicity bound for near-group fusion categories of type\\ $\mathbb Z/4\mathbb Z+k$}

\author{Yue Meng and Zhiqiang Yu}
\date{}
\maketitle

%\subjclass[2020]{18M20, 18M15}
%\keywords{Near-group fusion category, Drinfeld center, central induction, Frobenius--Schur indicator}

\begin{abstract}
Let $\cC$ be a near-group fusion category of type $\mathbb Z/4\mathbb Z+k$. We show  $k\in\{0,3,4\}$, which then gives a complete characterization of near-group fusion categories determined by $\mathbb Z/4\mathbb Z$.

\end{abstract}


%{\small\quad\textbf{Mathematics Subject Classification (2020)} 18M20}

\section{Introduction}
The fusion rules of a near-group fusion category are determined by a finite group $G$ and a nonnegative integer $k$, with the additional simple object $\rho$ satisfying $\rho\otimes\rho=\bigoplus_{g\in G}g\oplus k\rho$. We denote this fusion  category by $G+k$. For every finite abelian group $G$, near-group fusion categories of  type $G+0$ are exactly Tambara--Yamagami categories \cite{TambaraYamagami}. A natural question is therefore to determine which positive values of $k$ admit a near-group fusion category of type $G+k$. Ostrik proved that, for $G=\{e\}$ and $G=\mathbb Z/2\mathbb Z$, the admissible levels are precisely $k\in\{0,1\}$ and $k\in\{0,1,2\}$, respectively \cite{Ostrik2003,Ostrik2015}. For $G=\mathbb Z/3\mathbb Z$, Larson showed that the admissible levels are exactly $k\in\{0,2,3,6\}$ \cite{Larson}. More recently, Schopieray classified the admissible levels for all nontrivial elementary abelian $2$-groups \cite[Theorem~1.1]{Schopieray2023}: a near-group fusion category of type $G+k$ exists if and only if $k=0$, or $G=\mathbb Z/2\mathbb Z$ and $k\in\{1,2\}$, or $G=\mathbb Z/2\mathbb Z \times \mathbb Z/2\mathbb Z$ and $k=4$.

Let $\cC$ be a near-group fusion category of type $\ZZ/4\ZZ+k$, 
$ \mathcal{O}(\cC)=\{\one,g,g^2,g^3,\rho\}$
and fusion rules are defined as
\begin{equation}\label{eq:fusion}
 g^4=\one,\quad g^j\rho=\rho g^j=\rho,\quad
 \rho^2=\sum_{j=0}^3g^j+k\rho,
 \quad k\in\NN.
\end{equation}
The main results of this paper are as follows:

\begin{theorem}[Theorem \ref{thm:main}]
Let $\cC$ be a fusion category with fusion rules
\eqref{eq:fusion}. Then
$
 k\in\{0,3,4\}.
$
\end{theorem}

For $k=0$, these are the
Tambara--Yamagami categories $\operatorname{TY}(\ZZ/4\ZZ,\chi,\tau)$,
where $\chi$ is a nondegenerate symmetric bicharacter on $\ZZ/4\ZZ$
and $\tau=\pm\frac12$ \cite{TambaraYamagami}. For $k=3$, a
categorification is provided by
$\operatorname{Rep}(\ZZ/5\ZZ\rtimes \ZZ/4\ZZ)$ and
$\ZZ/4\ZZ$ acts faithfully on $\ZZ/5\ZZ$. For $k=4$, Izumi's construction
gives unitary near-group categories of type $\ZZ/4\ZZ+4$,  two inequivalent
unitary categorifications are recorded in
\cite[Table~2]{EvansGannon} and are complex conjugate to one
another \cite[Example~12.12]{IzumiNearGroup}. Thus the possible values of $k$ in the pseudo-unitary setting
are exactly $0$, $3$, and $4$.

The method follows Larson's use of central induction, twist traces,
and sums of roots of unity \cite{Larson}. We also use the dimension
row formula recorded in \cite[Lemma 2.1]{Schopieray}. Its application
here is an integrality obstruction for a fusion matrix. The numerical
inequality in \cite[Theorem 5.1]{Schopieray} by itself does not exclude
$k=8$; this distinction is explained in Section~\ref{levelbound}.

\section{Central induction and preliminary reductions}\label{section2}

\subsection{Standard identities}
%Every near-group fusion category is Galois conjugate to a
%pseudo-unitary fusion category; see [10, Appendix A.1].
%Since Galois conjugation preserves the fusion rules, it suffices
%to prove the restriction on $k$ in the pseudo-unitary case.
%Hence, we equip $\mathcal C$ with its canonical spherical
%structure and $\mathcal Z(\mathcal C)$ with the induced spherical
%structure. 

For objects $U,V$ in a fusion category $\cC$, let
$[U,V]:=\dim\Hom(U,V)$. We recall the following standard facts about a spherical fusion category and its Drinfeld center $\cZ(\cC)$; see
\cite[Section 2]{Ostrik2015}, \cite{NgSchauenburg}.

 Let $F:\cZ(\cC)\to\cC$ be the forgetful functor, and let $I$ be its right adjoint. By \cite[Proposition 9.2.2]{EGNO}, adjunction and the induction formula give
\begin{equation}\label{eq:FI}
 F(I(V))=\bigoplus_{X\in\mathcal{O}(\cC)}X\otimes V\otimes X^*,
 \quad
 [I(U),I(V)]=[U,F(I(V))].
\end{equation}
Let $D:=\dim(\cC)$. Notice that the scalar twists $\theta_X$ of  simple objects  $X\in\cO(\cZ(\cC))$ are roots
of unity, let $\nu_n(V)$ be the $n$-th Frobenius-Schur indicator of object $V\in\cC$, then  
\cite[Theorem~4.1]{NgSchauenburg} states that
\begin{equation}\label{eq:indicator-trace}
\Tr(\theta_{I(V)}^n)=D\cdot\nu_n(V)\quad(n\geq1),
\end{equation}
Let $g$ be an invertible object of order $n$.
Then the vector space $\operatorname{Hom}(\mathbf1,g^{\otimes n})$ is one-dimensional.
Since $\nu_n(g)$ is the trace of an endomorphism of finite order
on this space, it is a root of unity. Hence $|\nu_n(g)|=1$.
Since $\dim I(g)=D$ and   dimensions of an objects are positive,
Equation 
\eqref{eq:indicator-trace} forces
$\theta_X^n=\nu_n(g)$ for every simple summand $X$ of $I(g)$.

Let $\cC$ be a fusion category with Grothendieck ring $K_0(\cC)$. Then the algebra $K_0(\cC)\otimes_{\ZZ}\CC$ is semisimple. Given an irreducible representation $\chi$ of $K_0(\cC)\otimes_\ZZ\mathbb{C}$, there is a positive algebraic integer  $f_\chi$, called the formal codegree of $\cC$, satisfying $\sum_{\chi}\frac{\dim\chi}{f_\chi}=1$ by \cite[proposition 2.10]{Ostrik2015}. For example, $\mathrm{FPdim}(\cC)$ and $\dim(\cC)$ are formal codegrees of $\cC$.
 Moreover, \cite[Theorem~2.13]{Ostrik2015} gives the 
decomposition 
\begin{equation}\label{eq:codegrees}
 I(\one)=\bigoplus_\chi A_\chi,
 \quad
 \dim A_\chi=\frac{D}{f_\chi},
 \quad
 \theta_{A_\chi}=1.
\end{equation}

\subsection{Computation of the bi-adjoint functor and twists}

Let $\cC$ be a near-group fusion category of type $\ZZ/4\ZZ+k$, and let $d=\dim(\rho)$. Then $d$ is the positive root of $x^2-kx-4$, i.e.
$
 d=\frac{k+\sqrt{k^2+16}}2$. 
Set $D=\dim(\mathcal{C})=4+d^2=8+kd$ and its Galois conjugate 
$D'=8+k(k-d)$.

Since the fusion ring is commutative, its formal codegrees
are the eigenvalues of multiplication by
$
 \sum_{X\in\mathcal{O}(\cC)}XX^*
 =4\cdot\mathbf1+\rho^2
$
on $K_0(\mathcal{C})\otimes_{\mathbb{Z}}\mathbb{C}$
\cite[Remark~2.11]{Ostrik2015}.
With respect to the ordered basis
$(\mathbf1,g,g^2,g^3,\rho)$, this operator has matrix
\begin{equation*}
 R=
 \begin{pmatrix}
 5&1&1&1&k\\
 1&5&1&1&k\\
 1&1&5&1&k\\
 1&1&1&5&k\\
 k&k&k&k&8+k^2
 \end{pmatrix}.
\end{equation*}
Then the formal codegrees of $\cC$ are precisely the roots of the characteristic
polynomial of $R$, which are
$
 D, D', 4, 4, 4.
$
In particular, Equation \eqref{eq:codegrees} means
\begin{equation}\label{eq:I1-dimension-bound}
 \dim X\leq\frac {D}{4},\quad \text{for all simple objects }X\subset I(\one).
\end{equation}
If $d$ is rational, it is an integer, then 
$k=0,3$.  
Assume  $d$ is irrational below, then it follows from \cite[Theorem~A.6(i)]{Ostrik2015} that $k=4t$.

\begin{proposition}\label{prop:automatic-MF}
Let $\cC$ be a near-group fusion category of type $\ZZ/4\ZZ+k$, then each $I(g^j)$ has exactly five distinct simple summands, for $0\leq j\leq 3$.
\end{proposition}

\begin{proof}
For $I(\one)$, it is proved by   Equation \eqref{eq:codegrees}.
For any $h\in\{g,g^2,g^3\}$, the induction formula gives
\begin{equation}\label{eq:FI-invertible-general}
F(I(h))=4h\oplus(\oplus_{j=0}^3g^j)\oplus  k\rho.
\end{equation}
 Write
$I(h)=\oplus_Z m_ZZ$ with $m_Z\in\NN$. By definition, we have $
 \sum_Z m_Z^2=5$.
If $I(h)$ were not multiplicity-free, the only possible decomposition would then be
$
 I(h)=X\oplus2Y,\ X\not\cong Y$, for some $X,Y\in\mathcal{O}(\mathcal{Z}(\mathcal{C})).
$
Since $[I(h),I(\one)]=1$ and $I(\one)$ is
multiplicity-free, we have
$
 1=[I(\one),X]+2[I(\one),Y].
$
Consequently, $X\subset I(\one)$ and $Y\not\subset I(\one)$.
In particular, $\theta_X=1$.

By \cite[Theorem 4.1]{NgSchauenburg}, we know 
$0=\nu_1(h)\dim(\cC)=\Tr(\theta_{I(h)})=\dim( X)+2\theta_Y\dim (Y).
$
Hence, $\dim(X)=|\dim (X)|=|2\dim Y|=2\dim(Y)$ as $\cC$ is pseudo-unitary, consequently,  $\theta_Y=-1$. Since
$\dim I(h)=D$, it follows that $
 D=\dim (X)+2\dim( Y)=2\dim (X)$.
Thus $\dim (X)=D/2$, contradicting Equation (\ref{eq:I1-dimension-bound}).
This contradiction proves multiplicity-freeness for every $h$.

Finally, Equation \eqref{eq:FI-invertible-general} shows that each induction has five
distinct simple summands,  any
two different inductions have exactly one summand in common.
\end{proof}
%%%%%%%%%%%%%%%%%%%%
%\section{Determining the four multiplicity-free inductions}
%\subsection{Computation of the bi-adjoint functor and twists}
For all $X\in\mathcal{C}$ and $Z\in\mathcal{Z}(\mathcal{C})$, the adjointness condition gives 
$
[I(X),Z]=[X,F(Z)]$. 
\begin{lemma}\label{lem:Ostrik2.13}
 There are pairwise non-isomorphic
simple objects
$
\mathbf1,A,B_{01},B_{02},B_{03}\in \mathcal{O}(\mathcal Z(\mathcal C))
$
such that
\[
I(\mathbf1)
\cong \mathbf1\oplus A\oplus B_{01}\oplus B_{02}\oplus B_{03}.
\]
Their dimensions and forgetful images are
\[
\begin{aligned}
\dim A=1+td,\quad
F(A)\cong\mathbf1\oplus t\rho,\quad
\dim B_{0j}=2+td,\quad
F(B_{0j})\cong\mathbf1\oplus g^j\oplus t\rho
\;(1\le j\le3).
\end{aligned}
\]
\end{lemma}

\begin{proof}
It follows from Equation \eqref{eq:codegrees} and Proposition \ref{prop:automatic-MF} that the decomposition $I(\one)=\one\oplus A\oplus C_1\oplus C_2\oplus C_3$ is multiplicity-free  and 
\[
\dim A=\frac{D}{D'}=\frac{d^2}{4}=1+td,
\qquad
\dim C_i=\frac{D}{4}=2+td.
\]
%Here we used $d^2=4+4td$ and $d(k-d)=-4$.
Notice that $d$ is irrational,  so $F(A)$ contains
exactly one invertible summand and $t$ copies of $\rho$, whereas
each $F(C_i)$ contains exactly two invertible summands and $t$
copies of $\rho$.
Meanwhile, Equation \eqref{eq:FI} implies 
$
F(I(\mathbf1))
\cong
4\cdot\mathbf1\oplus\left(\bigoplus_{j=0}^{3}g^j\right)
\oplus4t\rho.
$
Each of $g,g^2,g^3$ occurs exactly once in this decomposition.
Thus we may label the
objects $C_i$ as $B_{01},B_{02},B_{03}$ so that
$
F(B_{0j})\cong\mathbf1\oplus g^j\oplus t\rho        \
(1\le j\le3).
$
\end{proof}
By Proposition~\ref{prop:automatic-MF}, for any simple object $X\in\mathcal{O}(\cZ(\cC))$, every
invertible summand of $F(X)$ has multiplicity at most one. Let $k=4t$ with $t\geq 1$ as in Section~\ref{section2}.
Denote the  simple objects of $\mathcal Z(\mathcal C)$ by
$
 \mathbf1, A, B_{ij}\ (0\leq i<j\leq3),
 P_j,Q_j\ (1\leq j\leq3),
 L_\alpha\ (\alpha\in\mathcal J),
$
where $\mathcal J$ is a finite index set. Here  $B_{ij}$ denotes the unique common
simple summand of $I(g^i)$ and $I(g^j)$; the six objects
$B_{ij}$ are shown to be pairwise non-isomorphic in the
following proposition.

\begin{proposition}\label{prop:central-inductions}

There are integers $p,q,r,s\geq0$ and $\gamma_\alpha>0$, and
$\eta\in\{i,-i\}$, such that
$
 p+q=r+s=t.
$
The block above the twist row is the matrix $W$ of
the forgetful functor $F$, with entries $W_{V,Y}=[F(Y),V]$.
Thus each column records the decomposition of $F(Y)$,
by adjunction, each row records the decomposition of $I(V)$. In particular,
the first four rows determine the decompositions of
$I(\mathbf1),I(g),I(g^2),I(g^3)$.
The last row records the twists and is not part of $W$.
The final column is repeated for every $\alpha\in\mathcal J$,
and each $\theta_\alpha$ is a root of unity. Zero entries are omitted.
\begingroup
\setlength{\arraycolsep}{6pt}
\renewcommand{\arraystretch}{1.2}
\begin{equation*}\label{eq:central-restrictions-twists}
\begin{array}{c|ccccc|cccc|ccc|cc|c}
 &\mathbf1&A&B_{01}&B_{02}&B_{03}
 &B_{12}&B_{13}&P_1&Q_1
 &B_{23}&P_2&Q_2&P_3&Q_3&L_\alpha\\ \hline
 \mathbf1&1&1&1&1&1&&&&&&&&&&\\
 g&&&1&&&1&1&1&1&&&&&&\\
 g^2&&&&1&&1&&&&1&1&1&&&\\
 g^3&&&&&1&&1&&&1&&&1&1&\\ \hline
 \rho&&t&t&t&t&t&t&p&q&t&r&s&p&q&\gamma_\alpha\\ \hline\hline
 \theta&1&1&1&1&1&-1&\eta&-\eta&-\eta&-1&1&1&-\eta&-\eta&\theta_\alpha
\end{array}
\end{equation*}
\endgroup
\end{proposition}

\begin{proof}
By Proposition~\ref{prop:automatic-MF}, $I(g^j)$ are multiplicity-free and
each have five simple summands, and meet pairwise in exactly
one summand. 
The decomposition of $I(\mathbf1)$ is obtained in Lemma \ref{lem:Ostrik2.13}, and  $I(g^j)$ contains a simple summand whose twist is $1$.
The twist-power constraint from Section~\ref{section2} gives
\begin{equation*}
 \theta_Y^2=1\quad(Y\subset I(g^2)),\qquad
 \theta_Y^4=1\quad(Y\subset I(g)\text{ or }I(g^3)).
\end{equation*}
The three pairwise intersections among $I(g),I(g^2),I(g^3)$
are either three distinct simple objects or a single simple
object $T$ common to all three inductions.
In the latter case, the invertible summands of $F(T)$ are
exactly $g,g^2,g^3$.
Moreover, $I(g)$ consists of $B_{01}$, $T$, and three further
simple objects whose restrictions contain $g$ as their only
invertible summand.
Denote their twists by $\xi_1,\xi_2,\xi_3$.
Since $\theta_T^2=1$ and $d\notin\mathbb Q(i)$, comparison of
coefficients in $\operatorname{Tr}(\theta_{I(g)}^2)=0$ gives
$
 2+3+\xi_1^2+\xi_2^2+\xi_3^2=0,\
  \xi_\ell^2\in\{1,-1\},
$
which is impossible.
Thus the three intersections are distinct.
Denote the common summand of $I(g^i)$ and $I(g^j)$ by
$B_{ij}$ for $1\leq i<j\leq3$.
The invertible summands of $F(B_{ij})$ are exactly $g^i,g^j$.
Each $I(g^j)$ therefore has exactly two remaining summands,
denoted $P_j,Q_j$, whose restrictions contain $g^j$ as their
only invertible summand.
This proves the first four rows of the table.

Set $\varepsilon=\theta_{B_{12}}\in\{1,-1\}$,
$\eta=\theta_{B_{13}}$, $\xi=\theta_{P_1}$, and
$\xi'=\theta_{Q_1}$.
The coefficients of $1$ in the second and first twist traces
of $I(g)$ give, respectively,
$
 4+2\eta^2+\xi^2+(\xi')^2=0,\
 2+2\varepsilon+2\eta+\xi+\xi'=0.
$
The first equation forces
$\eta^2=\xi^2=(\xi')^2=-1$.
The real part of the second then gives $\varepsilon=-1$,
and its imaginary part gives $\xi=\xi'=-\eta$.
The same argument for $I(g^3)$ yields
$\theta_{B_{23}}=-1$ and $\theta_{P_3}=\theta_{Q_3}=-\eta$.
Finally, the coefficient of $1$ in the first twist trace of
$I(g^2)$ is
$
 2-2-2+\theta_{P_2}+\theta_{Q_2}=0,
$
so $\theta_{P_2}=\theta_{Q_2}=1$.

Let $c_{12},c_{13}$, $p,q$ be the multiplicities of $\rho$  in $F(B_{12}),F(B_{13})$,
 $F(P_1),F(Q_1)$, respectively.
The $\rho$ coefficient in $F(I(g))$ and the $d$-coefficients
in the second and first twist traces give
 \begin{align*}
 t+c_{12}+c_{13}+p+q&=4t,\\
 t+c_{12}-c_{13}-p-q&=0,\\
 t-c_{12}+\eta(c_{13}-p-q)&=0.
 \end{align*}
Hence, $c_{12}=c_{13}=p+q=t$.
Duality identifies $B_{23}=B_{12}^*$ and allows the labelling
$P_3=P_1^*$, $Q_3=Q_1^*$, proving the corresponding entries
in the $\rho$ row.
Writing $r,s$ for the $\rho$ multiplicities of $F(P_2),F(Q_2)$,
the $\rho$ coefficient of $F(I(g^2))$ gives $3t+r+s=4t$.
Thus $r+s=t$.

The restriction of the  remaining simple objects  $L_\alpha$ ($\alpha\in \mathcal{J}$) do not contain
  invertible summands,  so $F(L_\alpha)=\gamma_\alpha\rho$
for some $\gamma_\alpha\in\mathbb Z_{>0}$.
Its twist $\theta_\alpha$ is a root of unity.
This completes the table.
\end{proof}

Notice that $\rho$ is self-dual, so the second Frobenius-Schur indictor $\nu_2(\rho)=\pm1$ \cite{NgSchauenburg}. 
\begin{proposition}\label{prop:moments}
Let $\nu=\nu_2(\rho)\in\{1,-1\}$, $
 u=pq$, $v=rs$, $w=\sqrt{t^2+1}$,
 $N=\sum_\alpha\gamma_\alpha^2$,
$S_j=\sum_\alpha\gamma_\alpha^2\theta_\alpha^j$, $j=1,2$. Then  we have \begin{align}
 N&=8+6t^2+4u+2v,\label{eq:N}\\
 S_1&=-t^2+2v-2tw+\eta(t^2-4u),\label{eq:S1}\\
 S_2&=-4(t-\nu)w-4u+2v.\label{eq:S2}
\end{align}
In particular,
\begin{equation}\label{eq:N-upper}
 0\leq u,v\leq\frac{t^2}{4},\quad
 N\leq8+\frac{15}{2}t^2.
\end{equation}
\end{proposition}

\begin{proof}
By Proposition \ref{prop:central-inductions}, $
I(\rho)=tA\oplus_{i<j}tB_{i,j}
 \oplus p(P_1\oplus P_3)\oplus q(Q_1\oplus Q_3)
\oplus rP_2\oplus sQ_2\oplus_\alpha\gamma_\alpha L_\alpha$.
Equations \eqref{eq:FI} implies
$F (I(\rho))=4t(\oplus_{j=0}^3g^j)\oplus(8+16t^2)\rho$, which then gives
\begin{equation*}
8+16t^2=[F (I(\rho)),\rho]=[I(\rho),I(\rho)]=7t^2+2(p^2+q^2)+r^2+s^2+N.
\end{equation*}
Using $p^2+q^2=t^2-2u$ and $r^2+s^2=t^2-2v$, we obtain Equation \eqref{eq:N}.

Substituting the twists from Proposition \ref{prop:central-inductions} into the identities $\text{Tr}(\theta_{I(\rho)})$ and $\text{Tr}(\theta_{I(\rho)}^2)$ yields
 \begin{equation*}
 0=4t+d\bigl(3t^2-2v+\eta(4u-t^2)+S_1\bigr),\quad
 \nu D=8t+d\bigl(4t^2+4u-2v+S_2\bigr).
 \end{equation*}
After simplification, we prove Equations \eqref{eq:S1} and
\eqref{eq:S2}. The  Inequalities \eqref{eq:N-upper} follow from
$p+q=r+s=t$.
\end{proof}

\section{Bounds for sums of roots of unity}
\label{section3}

We establish the arithmetic estimates for the identities from 
Section \ref{section2}. 
We write $\varphi$ for Euler's function and $\mu$ for the
M\"obius function, $\varphi(1)=\mu(1)=1$, and let
$\zeta_m=\exp(2\pi i/m)$. For a prime $p$ and a positive integer $m$, let $v_p(m)$ denote
the exponent of $p$ in the prime factorization of $m$.
\subsection{Normalized Galois averages}
\label{subsec:normalized-averages}

Let $c>1$ be square-free and set $K=\QQ(\sqrt c)$. Its conductor is \(f=c\) if \(c\equiv1\pmod4\), and \(f=4c\) if \(c\equiv2,3\pmod4\).
In particular, $K\subseteq\QQ(\zeta_m)$ if and only if $f\mid m$.

Choose a cyclotomic field $E$ containing $K$ and all roots under
consideration. For $w\in E$, define its normalized relative trace by
\begin{equation*}
 T_K(w)=\frac{\Tr_{E/K}(w)}{[E:K]}
       =\frac1{|G|}\sum_{\sigma\in G}\sigma(w),
 \quad G=\Gal(E/K).
\end{equation*}
Obviously, $T_K(-)$ is $K$-linear
and satisfies $T_K(a)=a$ for $a\in K$. Moreover, for every $\tau\in G$, $
 \tau(T_K(w))
 =\frac1{|G|}\sum_{\sigma\in G}(\tau\sigma)(w)=T_K(w)$ and $
  T_K(\tau(w))=T_K(w)
$, so  $T_K(w)\in E^G=K$. 

The normalization makes $T_K(w)$ independent of the choice of
cyclotomic field containing $K$ and $w$.
Indeed, if $E\subseteq E'$ and $w\in E$, transitivity of the trace gives
$
 \frac{\Tr_{E'/K}(w)}{[E':K]}
 =\frac{\Tr_{E/K}([E':E]w)}{[E':E][E:K]}
 =\frac{\Tr_{E/K}(w)}{[E:K]}.
$
Thus we may compute $T_K(w)$ in any cyclotomic field containing $K$
and $w$.
Since $K\subset\mathbb R$, complex conjugation belongs to $G$.
Consequently
\begin{equation*}
 T_K(\operatorname{Re}w)
 =T_K\!\left(\frac{w+\overline w}{2}\right)
 =\frac{T_K(w)+T_K(\overline w)}2=T_K(w).
\end{equation*}

\iffalse
Write $T_K(\xi)=x_\xi+y_\xi\sqrt c$ with $x_\xi,y_\xi\in\QQ$,
and put $O=O_K(\xi)$. Choose a positive integer $M$ such that
$Mx_\xi,My_\xi\in\ZZ$. Repeating each element of $O$ exactly
$M$ times gives a sum of $M|O|$ roots equal to
$
 M\sum_{\eta\in O}\eta
 =M|O|x_\xi+M|O|y_\xi\sqrt c.
$
Both coefficients on the right are integers. Theorem~\ref{thm:larson-root-bound}
 gives $M|O|\geq M|O|\,|y_\xi|\varphi(2c)$, and hence
\begin{equation}\label{eq:normalized-coefficient}
 |y_\xi|\leq\frac1{\varphi(2c)}.
\end{equation}
\fi
%\paragraph{Rational averages and coprime factors.}
If $\beta$ is a primitive $r$th root of unity, then
 \begin{equation*}
 \frac{\Tr_{\QQ(\beta)/\QQ}(\beta)}{\varphi(r)}
 =\frac1{\varphi(r)}
   \sum_{\substack{1\leq j\leq r,(j,r)=1}}\beta^j
 =\frac{\mu(r)}{\varphi(r)}.
 \end{equation*}
Since $\beta$ is a primitive $r$th root of unity, the classical
evaluation of Ramanujan sums
\cite[Section~1]{RoblesRoy2017} gives
$
\sum_{\substack{1\le j\le r,(j,r)=1}}\beta^j=\mu(r).
$ For every positive integer $r$, we have $|\mu(r)|\leq1$ and
$\varphi(r)\geq1$. Hence
$
 \frac{\mu(r)}{\varphi(r)}
 \leq\frac{|\mu(r)|}{\varphi(r)}
 \leq1,
$
with equality if and only if $r=1$.

Let $\ell,r$ be positive integers such that $f\mid\ell$ and
$(\ell,r)=1$, and put
$K_1=\QQ(\zeta_{\ell})$ and $K_2=\QQ(\zeta_r)$. Let $K_1K_2$ denote the compositum of $K_1$ and $K_2$, namely the smallest
subfield of $\mathbb C$ containing both. By the Chinese remainder theorem,
$
 \operatorname{Gal}(K_1K_2/\mathbb Q)
 \cong
 \operatorname{Gal}(K_1/\mathbb Q)
 \times
 \operatorname{Gal}(K_2/\mathbb Q).
$ Since $K\subseteq K_1$, an automorphism of $K_1K_2$ fixes $K$
if and only if its restriction to $K_1$ fixes $K$.
Its restriction to $K_2$ can be any element of
$\operatorname{Gal}(K_2/\mathbb Q)$. Hence
$
 \operatorname{Gal}(K_1K_2/K)
 \cong
 \operatorname{Gal}(K_1/K)\times
 \operatorname{Gal}(K_2/\mathbb Q).
$
Writing these two factors as $G_1$ and $G_2$, respectively,
for $u\in K_1$ and $v\in K_2$,  we obtain the formula
\[
 \begin{aligned}
 T_K(uv)
 &=\frac1{|G_1||G_2|}
   \sum_{\sigma\in G_1}\sum_{\tau\in G_2}\sigma(u)\tau(v)\\
 &=\left(\frac1{|G_1|}\sum_{\sigma\in G_1}\sigma(u)\right)
   \left(\frac1{|G_2|}\sum_{\tau\in G_2}\tau(v)\right)
 =T_K(u)\frac{\Tr_{K_2/\QQ}(v)}{[K_2:\QQ]}.
 \end{aligned}
\]

We will use this factorization only for coprime cyclotomic factors.

For a root of unity $\xi$, write
$
T_K(\xi)=x_\xi+y_\xi\sqrt c,
\ x_\xi,y_\xi\in\mathbb Q,
$
and let $m=\operatorname{ord}(\xi)$. If $f\nmid m$, then $K\not\subseteq\QQ(\xi)$ and $K\cap\QQ(\xi)=\QQ$.
Restriction identifies $\Gal(K\QQ(\xi)/K)$ with
$\Gal(\QQ(\xi)/\QQ)$, so the $K$-orbit of $\xi$ contains all
primitive $m$th roots. Therefore
$
 T_K(\xi)=\frac{\mu(m)}{\varphi(m)}\in\QQ.
$

Next assume that $f\mid m$ and that $v_p(m)>v_p(f)$ for a prime
factor $p$ of $ f$. These assumptions imply $p^2\mid m$ and $f\mid m/p$.
Every prime divisor of $m$ divides $m/p$, so $a=1+m/p$ is relatively
prime to $m$. The automorphism
$
 \sigma(\xi)=\xi^a
$
of $\QQ(\xi)$ fixes $K$, because $a\equiv1\pmod f$.
Let $\omega=\xi^{m/p}$. Then
$\sigma(\xi)=\omega\xi$ and
$\sigma(\omega)=\omega$, since $p\mid m/p$.
It follows by induction that $\sigma^j(\xi^k)=\omega^{jk}\xi^k$.
For every positive integer $k$ with $p\nmid k$, trace invariance
gives
\begin{equation}\label{pT}
 pT_K(\xi^k)=T_K(p\xi^k)
 =\sum_{j=0}^{p-1}T_K(\sigma^j(\xi^k))
 =T_K\!\left(\xi^k\sum_{j=0}^{p-1}\omega^{jk}\right)=0.
\end{equation}
In particular, $T_K(\xi)=0$. Combining rationality with this cancellation proves the necessary conditions: if \(y_\xi\neq 0\), then \(f\mid m\) and \(v_p(m)=v_p(f)\) for every prime \(p\mid f\).
%These conditions are not asserted to be sufficient.

%\subsection{A root of unity bound for even \texorpdfstring{$c$}{c}}
We use the following  theorem  without proof.

\begin{theorem}[{\cite[Theorem 6.2]{Larson}}]
\label{thm:larson-root-bound}
Let $c>1$ be square-free and let $a,b\in\ZZ$. If roots of unity
$\xi_1,\ldots,\xi_n$ satisfy
$
 \sum_{j=1}^n\xi_j=a+b\sqrt c,
$
then $n\geq |b|\varphi(2c)$.
\end{theorem}
The following lemma is a direct consequence of \cite[Lemma 3.11]{DongPalcouxWang2026}.
\begin{lemma}
\label{lem:even-moments}
Let $c>1$ be even and square-free. Suppose that roots of unity
$\xi_1,\ldots,\xi_n$ satisfy
\[
 \operatorname{Re}\sum_{j=1}^n\xi_j=x+y\sqrt c,
 \quad \sum_{j=1}^n\xi_j^2=z+q\sqrt c,
 \quad x,y,z,q\in\ZZ.
\]
Then
$
 n\geq\varphi(2c)(|y|+|q|).
$
\end{lemma}

\iffalse
\begin{proof}
Write $c=2h$ with $h$ odd. The conductor is $f=8h$, so
$v_2(f)=3$. For a root $\xi$ of order $m$, write
$
 T_K(\xi)=x_\xi+y_\xi\sqrt c,\
  T_K(\xi^2)=z_\xi+q_\xi\sqrt c.
$
The necessary conditions above give: if \(y_\xi\neq 0\), then \(v_2(m)=3\).
The order of $\xi^2$ is $m/(m,2)$. If $q_\xi\ne0$, that order
must have $2$-adic valuation $3$. In this case \(m\) is even, and \(v_2(m/(m,2))=v_2(m)-1\), so if \(q_\xi\neq 0\), then \(v_2(m)=4\).
Thus the two coefficients cannot both be nonzero.
Applying Inequality \eqref{eq:normalized-coefficient} to $\xi$ and $\xi^2$ yields
$
 y_\xi q_\xi=0,\
  |y_\xi|+|q_\xi|\leq\frac1{\varphi(2c)}.
$

By linearity and the identity $T_K(w)=T_K(\operatorname{Re}w)$,
the assumptions imply
 \begin{equation*}
 \sum_{j=1}^nT_K(\xi_j)=x+y\sqrt c,
 \quad \sum_{j=1}^nT_K(\xi_j^2)=z+q\sqrt c.
 \end{equation*}
Comparing the $\sqrt c$ coefficients and using the triangle
inequality, we obtain
$
 |y|+|q|
\leq\sum_{j=1}^n(|y_{\xi_j}|+|q_{\xi_j}|)
 \leq\frac n{\varphi(2c)}.
$
\end{proof}\fi
\subsection{A joint estimate over \texorpdfstring{$\QQ(\sqrt5)$}{Q(sqrt(5))}}
\label{subsec:five-moments}
Throughout this subsection, put $K=\QQ(\sqrt5)$ and
$T=T_K$.  Let $\zeta=\zeta_5$. The extension $\QQ(\zeta)/K$ has degree $2$, with complex
conjugation as its nontrivial automorphism. Its two orbits of
primitive fifth roots are therefore
$\{\zeta,\zeta^{-1}\}$ and $\{\zeta^2,\zeta^{-2}\}$.
Squaring interchanges them, so for every primitive fifth root
$\beta$ there is a sign $\varepsilon\in\{1,-1\}$ such that $
T(\beta)=\frac{-1+\varepsilon\sqrt5}{4},
 T(\beta^2)=\frac{-1-\varepsilon\sqrt5}{4}$.
The two signs are linked by the choice of $\beta$.

\begin{lemma}\label{lem:five-moments}
Let $\xi_1,\ldots,\xi_n$ be roots of unity such that
 \begin{equation*}
 \operatorname{Re}\sum_{j=1}^n\xi_j=x+y\sqrt5,
 \quad \sum_{j=1}^n\xi_j^2=z+q\sqrt5,
 \quad x,y,z,q\in\ZZ.
 \end{equation*}
Then
$
 n\geq-x-y-4q.
$
\end{lemma}
\begin{proof}
For a root of unity $\xi$, write
$
 T(\xi):=x_\xi+y_\xi\sqrt5,\
  T(\xi^2):=z_\xi+q_\xi\sqrt5,\
  \psi(\xi):=-x_\xi-y_\xi-4q_\xi.
$
We first prove $\psi(\xi)\leq1$ by considering $s=\ord(\xi)$.

If $5\nmid s$, then neither $s$ nor $\ord(\xi^2)$ is divisible
by the conductor $5$. Both averages are rational, so
$y_\xi=q_\xi=0$. Moreover, $x_\xi$ is the real part of an average
of roots of unity, and hence $x_\xi\geq-1$.
Thus $\psi(\xi)=-x_\xi\leq1$.
If $25\mid s$, the cancellation argument in
Subsection~\ref{subsec:normalized-averages} and Equation \eqref{pT},
with $p=5$ and $k=1,2$, gives $T(\xi)=T(\xi^2)=0$.
Hence $\psi(\xi)=0$ in this case.

In the remaining case $5\mid s$ but $25\nmid s$.
Removing the power of $2$ and the single factor $5$ gives
$
 s=2^a5r, a=v_2(s)\geq0, (r,10)=1.
$
Suppose first that $a\geq3$. The integers
$b_1=1+s/2$ and $b_2=1+s/4$ are odd and are $1$ modulo every
odd prime divisor of $s$. Thus they are relatively prime to $s$.
They are also $1$ modulo $5$, so the automorphisms
$\sigma_i(\xi)=\xi^{b_i}$ fix $K$.
Since $\xi^{s/2}=-1$, we have
$
 \sigma_1(\xi)=-\xi,\
  \sigma_2(\xi^2)=-\xi^2.
$
Since $T(\sigma(w))=T(w)$ for every automorphism $\sigma$
fixing $K$, we obtain
\begin{equation*}
 T(\xi)=T(\sigma_1(\xi))=T(-\xi)=-T(\xi),
 \quad
 T(\xi^2)=T(\sigma_2(\xi^2))=T(-\xi^2)=-T(\xi^2).
\end{equation*}
Thus $T(\xi)=T(\xi^2)=0$.

It remains to treat $a\in\{0,1,2\}$. Let $\ell=2^a5$.
Since $(\ell,r)=1$, choose integers $b,e$ satisfying $br+e\ell=1$ and set
$
 \delta=\xi^{br}, \gamma=\xi^{e\ell}.
$
Then $\xi=\delta\gamma$, and the congruences
$br\equiv1\pmod \ell$ and $e\ell\equiv1\pmod r$ show that $\delta$
and $\gamma$ are primitive roots of orders $\ell$ and $r$,
respectively. For $r=1$, this means $\gamma=1$.
Because $r$ is odd, squaring permutes its primitive roots. Thus
$\gamma$ and $\gamma^2$ have the same rational average
$\lambda=\mu(r)/\varphi(r)$.
The coprime factorization from
Subsection~\ref{subsec:normalized-averages} gives
\begin{equation*}
 T(\xi)=\lambda T(\delta),
 \quad T(\xi^2)=\lambda T(\delta^2),
 \quad \psi(\xi)=\lambda \psi(\delta).
\end{equation*}
The last identity follows by comparing the rational and $\sqrt5$
coefficients in the first two identities. Furthermore,
$
 -\frac12\leq\lambda\leq1.
$
Indeed, $|\mu(r)|\leq1$. If $\mu(r)=-1$, then the odd integer
$r$ is greater than $1$. Hence $1$ and $r-1$ are distinct and
both coprime to $r$, so $\varphi(r)\geq2$.

We now compute $\psi(\delta)$ for the three possible orders.
For $\ell=5$, take $\delta=\beta$ in the fifth-root formulas above.
For $\ell=10$, write $\delta=-\beta$ with $\beta$ primitive of order
$5$; then $T(\delta)=-T(\beta)$ and $T(\delta^2)=T(\beta^2)$.
For $\ell=20$, put $\eta=\delta^5$ and $\beta=\delta^{-4}$.
These roots have orders $4$ and $5$, satisfy $\delta=\eta\beta$,
and have $\eta^2=-1$. The rational average of $\eta$ is zero,
so coprime factorization yields
\begin{equation*}
 T(\delta)=0,
 \quad T(\delta^2)=-T(\beta^2)
                 =\frac{1+\varepsilon\sqrt5}{4}.
\end{equation*}
Substitution into the definition of $\psi$ gives the following table:
\[
 \begin{array}{c|ccc|c}
 \ell&x_\delta&y_\delta&q_\delta&\psi(\delta)\\ \hline
 5&-\tfrac14&\tfrac\varepsilon4&-\tfrac\varepsilon4
     &\tfrac{1+3\varepsilon}{4}\\[3pt]
 10&\tfrac14&-\tfrac\varepsilon4&-\tfrac\varepsilon4
     &\tfrac{-1+5\varepsilon}{4}\\[3pt]
 20&0&0&\tfrac\varepsilon4&-\varepsilon
 \end{array}
 \qquad (\varepsilon\in\{1,-1\}).
\]
In each row the sign is determined by its associated fifth root
$\beta$. These values give $-3/2\leq \psi(\delta)\leq1$.
If $\lambda\geq0$, then
$\psi(\xi)=\lambda \psi(\delta)\leq\lambda\leq1$.
If $\lambda<0$, multiplication reverses the lower bound on
$\psi(\delta)$, and hence
$
 \psi(\xi)=\lambda \psi(\delta)
 \leq-\frac32\lambda\leq\frac34<1.
$
This proves $\psi(\xi)\leq1$ for every root of unity.

Finally, use one common cyclotomic field for all the given roots.
The identity $T(w)=T(\operatorname{Re}w)$ and the hypotheses give
 \begin{equation*}
 \sum_{j=1}^nT(\xi_j)
 =T\!\left(\operatorname{Re}\sum_{j=1}^n\xi_j\right)
  =T(x+y\sqrt5)=x+y\sqrt5,\quad
 \sum_{j=1}^nT(\xi_j^2)
 =T(z+q\sqrt5)=z+q\sqrt5.
 \end{equation*}
Since $1,\sqrt5$ are linearly independent over $\QQ$, comparison
of coefficients gives
$
 \sum_jx_{\xi_j}=x,\
 \sum_jy_{\xi_j}=y,\
 \sum_jq_{\xi_j}=q.
$
Therefore,
$
 -x-y-4q
 =\sum_{j=1}^n(-x_{\xi_j}-y_{\xi_j}-4q_{\xi_j})
 =\sum_{j=1}^n\psi(\xi_j)\leq n.
$
\end{proof}

%The following lemma determines the equality case of
%\cite[Theorem~6.2]{Larson} for the sum $-4-4\sqrt5$.
\begin{lemma}
\label{lem:five-equality}
If there are  roots of unity $\xi_j$ ($1\leq j\leq 16$) such that $
 \sum_{j=1}^{16}\xi_j=-4-4\sqrt5,
$
then their multiset consists of eight copies of $\zeta_5^2$ and
eight copies of $\zeta_5^3$. Consequently,
$
 \sum_{j=1}^{16}\xi_j^2=-4+4\sqrt5.
$
\end{lemma}

\begin{proof}
Let $T(\xi)=x_\xi+y_\xi\sqrt5$. We first identify all roots
for which $y_\xi=-1/4$.
If $5\nmid\ord(\xi)$, then $y_\xi=0$; if
$25\mid\ord(\xi)$, then $T(\xi)=0$.
Otherwise $\ord(\xi)=5r$ with $(r,5)=1$.
Choose integers $u,v$ such that $ur+5v=1$ and let
$\alpha=\xi^{ur}$,\ $\beta=\xi^{5v}$.
Then $\xi=\alpha\beta$, with $\alpha$ and $\beta$ primitive of
orders $5$ and $r$. The coprime factorization in
Subsection~\ref{subsec:normalized-averages} and the fifth-root
calculation in Subsection~\ref{subsec:five-moments} give
\begin{equation*}
 T(\xi)=T(\alpha)\frac{\mu(r)}{\varphi(r)}
       =\frac{\mu(r)}{\varphi(r)}
          \left(\frac{-1+\varepsilon\sqrt5}{4}\right),
 \quad \varepsilon\in\{1,-1\}.
\end{equation*}
Here $r$ need not be odd. In particular,
$
 y_\xi=\frac{\varepsilon\mu(r)}{4\varphi(r)}\geq-\frac14.
$
Equality requires $|\mu(r)|=\varphi(r)=1$, which is possible
only when $r=1$ or $r=2$.
For $r=1$, equality requires $\varepsilon=-1$, giving the roots
$\zeta_5^2,\zeta_5^3$ and rational coefficient $x_\xi=-1/4$.
For $r=2$, we have $\mu(2)=-1$, so equality requires
$\varepsilon=1$, giving the roots $-\zeta_5,-\zeta_5^4$ and
rational coefficient $x_\xi=1/4$.

Averaging the assumed identity and comparing coefficients yields
$
\sum_{j=1}^{16}y_{\xi_j}=-4,\
 \sum_{j=1}^{16}x_{\xi_j}=-4.
$
The first equality is a sum of sixteen terms each at least $-1/4$,
so every term equals $-1/4$.
Let $a$ and $b$ count the roots in
$\{\zeta_5^2,\zeta_5^3\}$ and
$\{-\zeta_5,-\zeta_5^4\}$, respectively, including multiplicities.
Then the number of roots and the rational coefficients give
$
 a+b=16,\ -\frac a4+\frac b4=-4.
$
Thus $a=16$, $b=0$, and every original root is either
$\zeta_5^2$ or $\zeta_5^3$.
Let their multiplicities be $a_2,a_3$. The original sum is real,
and conjugation exchanges these two roots. Subtracting the
conjugate sum gives
$
 (a_2-a_3)(\zeta_5^2-\zeta_5^3)=0.
$
Since the roots are distinct, $a_2=a_3=8$.
Finally,
$
 \sum_{j=1}^{16}\xi_j^2
 =8(\zeta_5^4+\zeta_5^6)
 =8(\zeta_5+\zeta_5^{-1})=-4+4\sqrt5.
$
\end{proof}

\section{Proof of the main theorem}\label{levelbound}

\begin{proposition}\label{prop:t-bound}
Let $\cC$ be a near-group fusion category of type $\ZZ/4\ZZ+k$, $k=4t$. Then  $t\leq2$.
\end{proposition}

\begin{proof}
Write $t^2+1=m^2c$, where $m\in\mathbb Z_{\geq1}$ and $c$ is
square-free. For every odd prime divisor $p$ of $c$, we have
$t^2\equiv-1\pmod p$ and $p\nmid t$. Fermat's little theorem gives
$
 1\equiv t^{p-1}\equiv(-1)^{(p-1)/2}\pmod p.
$
Thus $(p-1)/2$ is even, so $p\equiv1\pmod4$.
If $c$ is odd, it is a product of such primes, and hence
$c\equiv1\pmod4$.

By Proposition~\ref{prop:moments}, we have
\begin{equation}\label{eq:real-moments}
 \operatorname{Re}S_1=-t^2+2v-2tm\sqrt c,\quad
 S_2=-4u+2v-4(t-\nu)m\sqrt c.
\end{equation}
%In the following, we show $t\leq 2$. 

\smallskip\noindent\textit{Case 1: $c$ is odd and $c\neq5$.}
Every prime divisor of $c$ is congruent to $1$ modulo $4$.
Since $c>1$ is square-free and $c\ne5$, it has a prime divisor
$p_0\ne5$, and hence $p_0\ge13$.
Moreover, $(p-1)/\sqrt p>1$ for every prime divisor $p$ of $c$.
Therefore,
$
\frac{\varphi(2c)}{\sqrt c}
=\prod_{p\mid c}\frac{p-1}{\sqrt p}
\ge \frac{p_0-1}{\sqrt{p_0}}
\ge \frac{12}{\sqrt{13}}.
$
Note $\nu\leq1$,   applying Theorem~\ref{thm:larson-root-bound} to $S_2$, we get
\begin{equation*}
 N\geq4(t-1)m\varphi(2c)
 \geq\frac{48}{\sqrt{13}}(t-1)\sqrt{t^2+1}
 >13t(t-1).
\end{equation*}
An elementary calculation shows $13t(t-1)>8+15t^2/2$ for $t\geq 3$, which contradicts Inequality \eqref{eq:N-upper}.

\smallskip\noindent\textit{Case 2: $c$ is even and $c\neq2$.}
Write $c=2h$ with $h>1$. Then the smallest possible odd prime divisor of $h$ is $5$, so
$
 \frac{\varphi(2c)}{\sqrt c}
 =\sqrt2\prod_{p\mid h}\frac{p-1}{\sqrt p}
 \geq\frac8{\sqrt{10}}>\frac52.
$
Apply Lemma~\ref{lem:even-moments} to \eqref{eq:real-moments}, we obtain
\begin{equation}\label{eq:even-applied}
 N\geq (6t-4\nu)m\varphi(2c)
 \geq(6t-4)m\varphi(2c).
\end{equation}
Consequently, $N>\frac52(6t-4)t=15t^2-10t$.
For $t\geq3$, this is larger than $8+15t^2/2$, impossible.

\smallskip\noindent\textit{Case 3: $c=2$.}
Inequality \eqref{eq:even-applied} implies 
$
 N\geq\sqrt2(6t-4)\sqrt{t^2+1}
 >\sqrt2(6t-4)t.
$
For every $t\geq7$,
$
 \sqrt2(6t-4)t>8+\frac{15}{2}t^2.
$
Indeed, the difference is
$(6\sqrt2-15/2)t^2-4\sqrt2t-8$; at $t=7$ it is
$266\sqrt2-751/2>0$, and its derivative is positive for $t\geq7$.
Thus $t\leq6$. None of $t=3,4,5,6$ has $t^2+1$ equal to twice
a square.

\smallskip\noindent\textit{Case 4: $c=5$.}
Lemma~\ref{lem:five-moments} gives
$
 N\geq t^2-2v+(18t-16\nu)m.
$
Using $\nu\leq1$, the formula for $N$, and the bounds on $u,v$,
we obtain
\begin{equation}\label{eq:short-budget}
 (18t-16)m
 \leq8+5t^2+4u+4v
 \leq8+7t^2
 =35m^2+1.
\end{equation}
If $m\geq2$, then $t^2=5m^2-1>4m^2$. Since $t,m$ are integers,
this implies $t\geq2m+1$. But then
\begin{equation*}
 (18t-16)m
 \geq\bigl(18(2m+1)-16\bigr)m
 =36m^2+2m
 >35m^2+1,
\end{equation*}
contradicting  Inequality \eqref{eq:short-budget}. Thus $m=1$ and $t=2$.

In summary, we conclude that $t\leq2$.
\end{proof}

We show that $k=8$ can be excluded in the following proposition. 
\begin{proposition}\label{prop:exclude-eight}
The case $k=8$ is impossible.
\end{proposition}

\begin{proof}
Suppose $k=8$, and retain the notations of Section~\ref{section2}.
Set
$
 d=4+2\sqrt5,\ a=1+2d,\ b=1+d,\ D=8b.
$
Then $d^2=4a$, $b^2=5a$, $D^2=320a$, and $a+a^{-1}=18$.
Also $u=pq,\ v=rs\in\{0,1\}$, since $p+q=r+s=2$.
Let $s$ denote the normalized $S$-matrix of the center $\mathcal{Z}(\cC)$, with
$s_{\mathbf1,X}=(\dim X)/D$.

\medskip
\noindent\textit{Case 1: $v=1$.}
If $u=0$, there is a central invertible object $Z$ with
$F(Z)= g$, and hence $F(Z^2)= g^2$.
Since $Z^2$ is simple, the restriction table forces it to be
either $P_2$ or $Q_2$. However, $r=s=1$ gives
$F(P_2)= F(Q_2)= g^2\oplus\rho$, a contradiction.
Thus $u=v=1$.

For the nontrivial automorphism $\sigma$ of $\mathbb Q(\sqrt5)$,
the dimension-row formula \cite[Lemma~2.1]{Schopieray} gives
$s_{A,X}=\sigma(\dim X/D)$.
Apart from $\mathbf1,A$,  the dimensions of simple objects of $\mathcal{Z}(\cC)$ are $b$, $2b$,
or $\gamma_\alpha d$. Since
$\sigma(b/D)=b/D$ and $\sigma(d/D)=-d/D$, it follows that
$
 s_{A,A}=\frac1D,\
 s_{A,X}^{\,2}=s_{\mathbf1,X}^{\,2}\
 (X\notin\{\mathbf1,A\}).
$
Let $P$  be a simple object of  dimension $b$. Then
$s_{P,\mathbf1}=s_{P,A}=1/8$.
The Verlinde formula and orthogonality of the $s_{\mathbf1,-}$ and
$s_{P,-}$   give $
 N_{A,A}^{P}
 =\sum_X
   \frac{s_{A,X}^{\,2}-s_{\mathbf1,X}^{\,2}}{s_{\mathbf1,X}}
   \overline{s_{P,X}}
 =\frac{a^2-1+a^{-1}-a}{8D}
   =\frac{16(a+1)}{8D}=\frac12$,
where the last equalities use $a+a^{-1}=18$ and $D=4(a+1)$, it is a contradiction.

\medskip
\noindent\textit{Case 2: $v=0$.}
Relabel the simple objects so that $r=0$, and let $h=P_2$. Then
$F(h)=g^2$, $\theta_h=1$, and $h^2=\mathbf1$, since the unit is
the only central simple object restricting to $\mathbf1$.
For each simple object $X$, define $\varepsilon_X$ by
$
 c_{X,h}\circ c_{h,X}
 =\varepsilon_X\,\operatorname{id}_{h\otimes X},
$
where $c$ denotes the braiding. The hexagon and balancing equation give
$
 \varepsilon_X\in\{1,-1\},\
 \theta_{h\otimes X}=\varepsilon_X\theta_X,\
 s_{h,X}=\varepsilon_X\frac{\dim X}{D}.
$
Since $h\ne\mathbf1$, orthogonality of  $s_{h,-}$ and $s_{\one,-} $ implies
$
 \sum_{\varepsilon_X=-1}\dim(X)^2=\frac{D^2}{2}=160a.
$

Tensoring by $h$ maps each invertible summand $g^j$ of
$F(X)$ to $g^{j+2}$, with exponents taken modulo $4$.
The table of Proposition \ref{prop:central-inductions} shows that it preserves the twists of objects whose
restrictions contain exactly one invertible summand.
It also fixes $B_{02}$ and $B_{13}$, since each is uniquely
determined by the invertible summands of its restriction.
All these objects have positive sign.
Tensoring by $h$ exchanges $B_{01}$ with $B_{23}$ and
$B_{03}$ with $B_{12}$, and the two objects in each pair
have opposite twists.
Thus, among the previously determined objects, precisely
these four have negative sign, and each has dimension $2b$.

Let
\begin{equation*}
N_\pm=\sum_{\varepsilon_{L_\alpha}=\pm1}\gamma_\alpha^2,\qquad
 S_{j,\pm}=\sum_{\varepsilon_{L_\alpha}=\pm1}
                 \gamma_\alpha^2\theta_\alpha^j
 \quad(j=1,2).
\end{equation*}
The preceding dimension count and the identity $N=32+4u$
of Equation \eqref{eq:N} give
\begin{equation}\label{eq:eight-sign-budgets}
 N_-=\frac{160a-4(2b)^2}{d^2}=20,\quad
 N_+=12+4u\leq16.
\end{equation}
Tensoring by $h$ pairs the negative-sign $L_\alpha$ with equal
multiplicities $\gamma_\alpha$ and opposite twists. Indeed,
$g^2\rho=\rho$, and a fixed point would have
$\theta_{L_\alpha}=-\theta_{L_\alpha}$.
Hence $S_{1,-}=0$, and the formula for $S_1$ gives
\begin{equation*}
 S_{1,+}=S_1=-4-4\sqrt5+4\eta(1-u).
\end{equation*}
Applying Theorem~\ref{thm:larson-root-bound} to
$S_{1,+}+\overline{S_{1,+}}=-8-8\sqrt5$ gives
$2N_+\geq8\varphi(10)=32$. Thus $N_+=16$ and $u=1$.
Lemma~\ref{lem:five-equality} now applies to
$S_{1,+}=-4-4\sqrt5$ and gives $S_{2,+}=-4+4\sqrt5$.
Using the formula for $S_2$ in Proposition \ref{prop:moments} with $u=1$ and $v=0$,
and subtracting the preceding expression for $S_{2,+}$, we obtain
$
 S_{2,-}=S_2-S_{2,+}=-4(3-\nu)\sqrt5.
$
This is a sum of $N_-=20$ roots of unity, counted with
multiplicity. Theorem~\ref{thm:larson-root-bound} therefore implies
$
 20\geq4(3-\nu)\varphi(10)=16(3-\nu)\geq32,
$
a contradiction.
\end{proof}
Proposition \ref{prop:t-bound} and 
Proposition~\ref{prop:exclude-eight} together imply the following theorem.
\begin{theorem}\label{thm:main}
There is a near-group fusion category of type $\ZZ/4\ZZ+k$ if and only if $k\in\{0,3,4\}$.
\end{theorem}

\subsection*{Acknowledgements}
 Z. Yu  is supported
by NSFC (no.12571041) and Qinglan Project of Yangzhou University.
The authors acknowledge the use of AI as an auxiliary tool during the preparation of this manuscript, including for mathematical discussion and editorial assistance. All
mathematical content was independently verified by the authors, who take full responsibility for
the final manuscript.

\author{Yue Meng\\ \thanks{Email:\,DX120250050@stu.yzu.edu.cn}
\\{\small School  of Mathematics,  Yangzhou University,
Yangzhou 225002, China}
}\\

\noindent \author{Zhiqiang Yu\\ \thanks{Email:\,zhiqyumath@yzu.edu.cn}
\\{\small School  of Mathematics,  Yangzhou University,
Yangzhou 225002, China}
}
\end{document}